\documentclass[11pt]{article}
\usepackage{amsmath}
\usepackage{amsthm}
\usepackage{amssymb}
\usepackage{amscd}
\usepackage[dvipdfm]{hyperref}

\newcommand{\C}{\mathbb{C}}

\newcommand{\N}{\mathbb{N}}
\theoremstyle{plain}
\newtheorem{theorem}{Theorem}[section]
\newtheorem{proposition}[theorem]{Proposition}
\newtheorem{lemma}[theorem]{Lemma}
\newtheorem{Corollary}[theorem]{Corollary}

\theoremstyle{definition}

\makeatletter
\@addtoreset{equation}{section}
\makeatother
\newcounter{constantLABEL}
\newcommand{\cref}[1]{C_{\ref{#1}}}

\newcounter{constantslabel}
\begin{document}


\title{
A weak compactness theorem of the Donaldson--Thomas instantons 
on compact K\"{a}hler threefolds
} 
\author{Y. Tanaka}
\date{}


\maketitle


\begin{abstract}
In \cite{tanaka1}, 
we introduced a gauge-theoretic equation on symplectic
 6-manifolds, which is a version of the Hermitian--Einstein equation
 perturbed by Higgs fields, and called it a  {\it Donaldson--Thomas equation}, 
to analytically approach the Donaldson--Thomas invariants. 
In this article, we consider the equation on compact 
K\"{a}hler threefolds, and study 
 some of the analytic properties of solutions to them, 
using analytic methods in higher-dimensional Yang--Mills theory 
developed by Nakajima
 \cite{Nakajima87}, \cite{Nakajima88} and  
Tian \cite{Tian00}  
with some additional 
arguments concerning an extra nonlinear term coming from the
 Higgs fields.   
We prove that 
a sequence of solutions to 
the Donaldson--Thomas equation 
of a unitary vector bundle over a compact K\"{a}hler threefold 
has a converging subsequence outside a closed subset whose
real two-dimensional Hausdorff measure is finite, 
provided that the $L^2$-norms of the Higgs fields are uniformly bounded. 
We also prove an $n/2$-compactness theorem of solutions to the equations on compact
 K\"{a}hler threefolds.  
\end{abstract}

\markboth{}
{A weak compactness theorem of the Donaldson--Thomas instantons}


\section{Introduction}

The Donaldson--Thomas invariant is a 
 deformation invariant of Calabi--Yau threefolds,  
which was constructed by Thomas \cite{Thomas00}   
from the moduli space of (semi-)stable sheaves 
by using algebraic geometry techniques. 
There are further generalization of this 
by Joyce and Song \cite{JS} and Kontsevich and Soibelman \cite{KS},
\cite{KS2}. 
These fit into programmes by Donaldson and Thomas \cite{DoTh98} and 
 Donaldson and Segal \cite{DS0902}, 
and many outcomes of them were made in both Mathematics and Physics.

We approach these invariants by using an analysis, aiming at revealing 
more symmetry and structures in the theory of the Donaldson--Thomas invariants. 
In \cite{tanaka1}, 
we introduced perturbed Hermitian-Einstein equations on symplectic
 6-manifolds, which we called {\it Donaldson--Thomas equations}, with a
 solution called {\it Donaldson--Thomas instanton},  
to analytically approach the Donaldson--Thomas invariants,  
studied the infinitesimal deformation 
and the Kuranishi model of the moduli space of the Donaldson--Thomas
instantons,  
and described the moduli space 
as a symplectic quotient by using a moment map for the action of gauge
group.  
We also introduced a stability condition which ought to produce a
Hitchin--Kobayashi-type correspondence for 
the Donaldson--Thomas instanton on K\"{a}hler threefolds.

In this article, 
to make the
argument simple in some sense, 
we assume that the underlying manifold is a compact K\"{a}hler threefold,  
and look into the analytic aspect of the Donaldson--Thomas (D--T) instantons, especially,  
bubbling phenomena of them  at the initial phase.  
Bubbling phenomena of Yang--Mills fields  
were first studied by Uhlenbeck \cite{Uh1},
\cite{Uh2} (see also \cite{KW}), and later by 
Nakajima \cite{Nakajima87}, \cite{Nakajima88}.  
Tian \cite{Tian00} further analysed them  
by using geometric measure theoretic methods  
developed by Lin \cite{Lin99}.  
We use these methods  
with some additional 
arguments concerning an extra nonlinear term coming from the
 Higgs fields to analyse the Donaldson--Thomas instantons on
compact K\"{a}hler threefolds.

\paragraph{The equations on compact K\"{a}hler threefolds.}

Let $Z$ be a compact K\"{a}hler threefold with K\"{a}hler form
$\omega$, and let $E$ be a unitary vector bundle over $Z$ of rank $r$. 
A complex structure on $Z$ gives the splitting of 
the space of the complexified two forms as $
 \Lambda^2 \otimes \C = \Lambda^{1,1} \oplus \Lambda^{2,0} 
\oplus \Lambda^{0,2}$, 
and $\Lambda^{1,1}$ further decomposes into $ \C \langle \omega \rangle \oplus \Lambda_{0}^{1,1}$.

We consider the following equations for a connection $A$ of $E$, 
 and an $\mathfrak{u} (E)$-valued (0,3)-form $u$ on $Z$:  
\begin{gather}
F_{A}^{0,2} =0 , \quad   \bar{\partial}_{A}^{*} u = 0, 
\label{DT1}
\\
F_{A}^{1 ,1} \wedge \omega^{2} + [ u , \bar{u}] + i \frac{\lambda(E)}{3} Id_{E}
 \omega^3 = 0  ,
\label{DT2}
\end{gather}
where $\lambda (E)$ is a constant defined by 
$\lambda (E) := 
6 \pi ( c_1(E) \cdot [\omega]^{2} ) / r [\omega]^{3}$. 
We call these equations the {\it Donaldson--Thomas
equations}, and we call a solution $(A, u)$ to these equations a  
{\it Donaldson--Thomas
instanton} (or a {\it D--T instanton} for short).

We remark that, in the K\"{a}hler case, 
the Weitzenb\"{o}ck formula \eqref{weit} 
implies that the Higgs field $u$ is absent for $c_1 (Z) > 0$, 
and covariantly constant for $c_1 (Z) =0$. 
This property of the Higgs field is similar to the Hitchin pair
\cite{Hit}; so this article may 
virtually concern compact K\"{a}hler threefolds of general type. 
However, we expect that this might provide a model in some sense for attacking
problems on compact symplectic 6-manifolds with $c_1 =0$   
because of strong similarities between these two geometries.

As mentioned above, 
using the methods by  \cite{Nakajima87},  \cite{Nakajima88}, 
\cite{Tian00} 
with some additional 
arguments on the
 Higgs fields, 
we  prove   
a weak compactness theorem of the Donaldson--Thomas
instantons on compact K\"{a}hler threefolds. 
\begin{theorem}
Let $Z$ be a K\"{a}hler threefold, and let $E$ be a 
unitary vector bundle over
 $Z$. 
Let $\{ (A_n , u_n) \}$ be a sequence of D--T instantons of E. 
We assume that  $\int_{Z} |  u_{n} |^2 dV_{g}$ are uniformly bounded. 
Then there exists a subsequence $\{ (A_{n_{j}} , u_{n_{j}}) \}$ of $\{ (A_n , u_n) \}$,  
a closed subset $S$ of $Z$ whose real two-dimensional Hausdorff
 measure is finite, and a sequence of gauge transformations 
$\{ \sigma_{j} \}$over $Z \setminus S$
 such that 
$\{ \sigma_{j}^{*} (A_{n_{j}} , u_{n_{j}}) \}$ converges to a
 D--T instanton over 
$Z \setminus S$. 
\label{th:main}
\end{theorem}

In addition, 
following \cite{MR1926448}, \cite{MR2086545},  
we also prove an $n/2$-compactness theorem 
of the Donaldson--Thomas
instantons on compact K\"{a}hler threefolds.

\begin{theorem}
 Let $\{ (A_n , u_n) \}$ be a sequence of D--T instantons of
 a unitary vector bundle $E$ over a compact K\"{a}hler threefold $Z$ with 
$ \int_{Z} |F_{A_n}|^3 dV_{g} \leq C$ ,
where 
$C > 0$ is a uniform constant. 
We assume that $\int_{Z} | u_{n} |^2 dV_{g}$ 
are also uniformly bounded. 
Then there exists a sequence of 
gauge transformations $\{ \sigma_{j} \}$ and 
a subsequence 
$\{ (A_{n_{j}} , u_{n_{j}}) \}$ of 
$\{ (A_n , u_n) \}$ such that 
$\{ \sigma_{j}^{*} (A_{n_{j}} , u_{n_{j}}) \}$ converges to a smooth 
 D--T instanton of $E$ over $Z$. 
\label{th:conv}
\end{theorem}

\paragraph{On the assumption on the uniform bound on the $\bf L^2$-norm of
    $\bf u$.}

As in the case of the Hitchin pair \cite{Hit}, 
there is a circle action on the
moduli space of D--T instantons by $(A , u) \mapsto (A , e^{i\theta}
u)$ when the underlying manifold is a K\"{a}hler threefold. 
The action is Hamiltonian, and the moment map with respect 
to the action is given by $|| u
||_{L^2}^{2}$. 
For the Hitchin pair, Hausel \cite{Hau} introduced a compactification of
the moduli space of the Hitchin pairs by using the symplectic cut
developed by Lerman \cite{Ler}.  
We pursue an analogy of this to compactify the moduli space of D--T
instantons in the direction of the Higgs fields;  
so assuming $L^2$-bound for $u$ should fit into the context.

\vspace{0.5cm}
The organization of this article is as follows.  
In Section \ref{sec:mono}, we prove a monotonicity formula 
for the Donaldson--Thomas instantons on
compact K\"{a}hler threefolds. 
In Section \ref{sec:est}, we derive a bound on $u$, and prove an
$\varepsilon$-regularity theorem for the Donaldson--Thomas instantons on
compact K\"{a}hler threefolds. 
The proof of Theorem \ref{th:main} is given in Section \ref{sec:weak}. 
In Section \ref{sec:conv}, 
we prove the $n/2$-compactness theorem (Theorem \ref{th:conv}) of the
Donaldson--Thomas instantons on compact K\"{a}hler threefolds.

\paragraph{Notation.} Throughout this article, 
$C, C' $, and $C''$ are positive constants, but they can be different each
time they 
occur.

\paragraph{Acknowledgements.}

I would like to thank Ryoichi Kobayashi and Hiroshi Konno for 
valuable comments, and Mikio Furuta for valuable comments 
and many intensive discussions. I am also grateful to Katrin Wehrheim 
for wonderful encouragement and useful comments.   
Last but not least I would like to thank Dominic Joyce for many useful
comments, 
 and for pointing out
errors in an earlier version of this article.

\section{Monotonicity formula}
\label{sec:mono}

Let $Z$ be a compact K\"{a}hler threefold with K\"{a}hler metric $g$, 
and let $E$ be a unitary vector bundle
over $Z$ of rank $r$. 
We fix a Hermitian metric on $E$.  
Our starting point is the following identity: 
\begin{equation*}
\begin{split}
& \qquad \int_{Z} \left\{ 
| F_{A}^{0 ,2} + \bar{\partial}_{A}^{*} u |^2 + 
 \frac{1}{2} |F_{A}^{1,1}
  \wedge \omega^2 + [u , \bar{u} ] + i \frac{\lambda (E)}{3} Id_{E} \omega^3  
 |^2 \right\} dV_{g} \\
& \quad \qquad \qquad \qquad + (2 c_2(E)- c_1(E)^2) \cdot  [\omega] 
 +  \frac{ 3(c_1(E) \cdot [\omega]^{2})^2}{2 r [\omega]^2 }  \\
& \quad \quad \qquad  = \int_{Z} \left\{ 
 \frac{1}{2} \left( | F_{A} |^2  + | \bar{D}^{*}_{A} u |^2   
   +  | [u , \bar{u}]|^2  \right) \right. \\
& \, \, \, \quad \qquad \qquad \qquad \qquad \qquad  
  \left.  + \langle F_{A}^{1,1} \wedge \omega^2  + i \frac{\lambda (E)}{3} 
  Id_{E} \omega^3   , 
      [ u , \bar{u} ] \rangle 
\right\} dV_{g}  , \\
\end{split}
\end{equation*}
where $\bar{D}_{A}^{*} u = \bar{\partial}_{A}^{*} u + \partial_{A}^{*}
\bar{u}$. 
Here we used the Bianchi identity $\bar{\partial}_{A}
F_{A}^{0,2} =0$ to deduce $ \int_{Y} \langle F_{A}^{0,2} , 
 \bar{\partial}_{A}^{*} u \rangle dV_{g} = 0$. 
We put 
\begin{equation*}
\begin{split}
L{(A,u)} &= \int_{Z} \left\{ 
 \frac{1}{2} \left( | F_{A} |^2  + | \bar{D}^{*}_{A} u |^2   
   + | [u , \bar{u}]|^2  \right) \right.  \\
  & \qquad \qquad \qquad \qquad  
\left.  + \langle F_{A}^{1,1} \wedge \omega^2  
+ i \frac{\lambda (E)}{3}  Id_{E} \omega^3   , 
      [ u , \bar{u} ] \rangle \right\} dV_{g} . 
\end{split}
\end{equation*}
If $(A,u)$ is a D--T instanton, $L(A,u)$ becomes 
\begin{equation*}
\begin{split} 
L(A,u) &= \frac{1}{2} \int_{Z} \left\{ 
  | F_{A}^{\perp} |^2  - |[ u , \bar{u}] |^2  \right\} dV_{g}  \\
 &=(2 c_2(E)- c_1(E)^2) \cdot  [\omega] 
   + \frac{ 3(c_1(E) \cdot [\omega]^{2})^2}{2 r [\omega]^2 } , \\
\end{split}
\end{equation*}
where $F_{A}^{\perp}$ is the $\Lambda^{1,1}_{0}$-component of $F_{A}$. 
One can see also from this that the 
$L^2$-norm of the curvature $F_{A}$ is bounded by the
topological constant and $\int_{Z} | [ u , \bar{u}]|^2 dV_{g}$, 
if $(A,u)$ is a D--T instanton.  
Hence, from Proposition 3.1, it is bounded if $|| u ||_{L^2}$ is bounded.

In this section, we prove the following monotonicity formula for the
Donaldson--Thomas instantons on K\"{a}hler threefolds.   
\begin{proposition}
Let $(A,u)$ be a D--T instanton of 
a unitary vector bundle $E$ over a compact K\"{a}hler threefold $Z$.  
Then, for any $z \in Z$, there exists a  positive constant $r_z$ such that 
for any $0< \sigma < \rho< r_z$, the following holds: 
\begin{equation}
\begin{split}
&\frac{1}{\rho^2} e^{ a \rho^2} \int_{B_{\rho}(z)} 
 m(A,u) dV_g 
 -  \frac{1}{\sigma^2} e^{
a \sigma^2} \int_{B_{\sigma}(z)} 
 m (A,u) dV_g \\
& \quad \qquad \geq   \int_{\sigma}^{\rho} 8 \tau^{-3} e^{a \tau^2}
  \int_{B_{\tau} (z)}  | [ u, \bar{u}]|^2 dV_{g} d \tau \\ 
  &\qquad \qquad +  \int_{B_{\rho} (z) \setminus B_{\sigma} (z)} r^{-2} e^{a
  r^2}  \left\{  4 \left|  \frac{\partial }{\partial r}  
 \lfloor F_{A}^{\perp} \right|^2 
 - 12 \left| \frac{\partial }{\partial r}  
 \lfloor [u , \bar{u}] \right|^2 \right\}  dV_{g} , \\ 
\label{eq:monotone}
\end{split}
\end{equation}
where 
${m} (A, u) := |F_{A}^{\perp}|^2 -  | [u, \bar{u} ]|^2$, and 
$a$ is a constant which depends only on $Z$.
\label{th:monotone}
\end{proposition}

\begin{proof}
The proof goes along almost the same line as that of \cite[Th.~2.1.1,
 2.1.2]{Tian00}
(see also \cite[Th.~1]{Price}).  
Thus, we describe it rather sketchily.

First, as in \cite[pp.~208]{Tian00}, 
we consider a one-parameter family of 
diffeomorphisms $\{ \varphi_t \}_{|t| < \infty}$ 
of $Z$ with $\varphi_{0} = id_Z$. 
We fix a connection $A_{0}$, 
and denote by $D$ its covariant
derivative. 
For $(A,u) \in \mathcal{A}(E) \times \Omega^{0,3} (Z, \mathfrak{u}(E))$, 
where $\mathcal{A} (E)$ is the space of connections of $E$,  
we define 
a one-parameter family $\{ (A_t  , u_t )  \}$  
in the following way. 
Let $\tau_{t}^{0}$ be the parallel transport of $E$ associated to $A_{0}$
along the path $\varphi_{\bar{t}} (z)_{0\leq \bar{t} \leq t}$, where $z
\in Z$. 
We define a family of connections $A_{t}$ by defining its covariant
 derivative as 
$ D_{X}^{t} s := (\tau_{t}^{0})^{-1} 
\left( D_{d \varphi_{t} (X)} \left( \tau_{t}^{0} (s) 
\right)\right)$, 
where $X \in TZ$ and $s \in \Gamma (Z ,E)$. 
Then the curvature of $A_{t}$ is written as 
$F_{A_{t}} (X,Y) = 
(\tau_{t}^{0})^{-1} \cdot 
F_{A} (d \varphi_{t}(X) , d \varphi_{t} (Y) ) \cdot \tau_{t}^{0}$, 
where $X, Y \in TZ$. 
We also define  
$u_{t}$ by $\varphi_{t}^{*} u $.  
We now assume that $(A,u)$ is a D--T instanton. 
Then, the same computation as in \cite[pp.~209]{Tian00} yields   
\begin{equation*}
\begin{split}
0 &= 
\left. \frac{d}{dt}   L (A_{t} , u_{t}) \right|_{t=0}  
\\ &= - \frac{1}{2} \int_{Z} 
\left( |F_{A}^{\perp}|^2 
\text{div} X + 4 \sum \langle 
F_{A}^{\perp} ( [X , e_i] , e_{j}) , 
F_{A}^{\perp} ( e_i , e_j ) \rangle  \right) dV_{g}  \\
& \quad \quad  
- \frac{1}{2} \int_{Z} 
\left( \left| [u , \bar{u} ] \right|^2 \text{div} X \right.  \\
& \quad \quad \quad  
\left. + 12 \sum \langle  [u , \bar{u}] ( [X , e_i] , e_{i_2} , \dots ,
  e_{i_6} ) , [u , \bar{u}] ( e_{i_1} , \dots ,
  e_{i_6}) \rangle 
\right) dV_{g} , \\
\end{split}
\end{equation*}
where $\{ e_i\}$ is an orthonormal basis of $T_{z}Z$ at $z \in Z$, 
and $X = \left. \frac{\partial \varphi_{t}}{\partial t} \right|_{t=0}$. 
Furthermore, 
using $[X , e_i] = \nabla_{X} e_i - \nabla_{e_i} X$ and 
$g(\nabla_{X} e_i , e_k) = - g( \nabla_{X} e_k , e_{i})$, 
we get 
\begin{equation}
\begin{split}
&0  =  \int_{Z} 
\left( |F_{A}^{\perp}|^2 
\text{div} X  - 4 \sum \langle 
F_{A}^{\perp} ( \nabla_{e_i} X  , e_{j}) , 
F_{A}^{\perp} ( e_i , e_j ) \rangle  \right) dV_{g}  \\
& \quad  
- \int_{Z} \left( \left| [u , \bar{u} ] \right|^2 \text{div} X \right.  \\
& \quad \quad 
\left. - 12 \sum \langle  [u , \bar{u}] ( \nabla_{e_i} X , e_{i_2} , \dots ,
  e_{i_6} ) , [u , \bar{u}] ( e_{i_1} , \dots ,
  e_{i_6} ) \rangle 
\right) dV_{g} , \\
\end{split}
\label{fvar}
\end{equation}
This is a version of 
the first variation formula for the Donaldson--Thomas instantons
 on a compact K\"{a}hler threefold, from which we deduce, for example, the
 monotonicity formula.

Next, for $z \in Z$, we take a positive number $r_{z}$ so that 
the following holds:  
there are normal coordinates $\zeta= (\zeta_{1} ,\dots ,\zeta_{6})$ 
with $z = (0, \dots , 0)$ in the geodesic
 ball $B_{r_{z}} (z)$ of $Z$ with respect to the metric $g$, and    
\begin{gather*}
 |g_{ij} - \delta_{ij} | \leq c(z) ( |\zeta_1|^2 + \cdots + |\zeta_6|^2 )
 , \quad  
 | d g_{ij} | \leq c(z) \sqrt{|\zeta_1|^2 + \cdots + |\zeta_6|^2 }   
\end{gather*}
for some constant $c(z) >0$ which depends only on $r_{z}$ and
 the curvature of $g$. 
We then denote by $r=r(\zeta)$ the distance function from $z$, and by $\phi$ a
 positive function on the unit sphere $S^5$, and 
define a cut-off vector field $X$ by 
$ X ( \zeta ) := \xi (r) \phi \left( \zeta / r \right) 
r \frac{\partial}{\partial r}$, 
where $\xi$ is a smooth function with compact support in $B_{r_{z}} (z) $. 
We now take an orthonormal basis $\{ \frac{\partial}{\partial r},  
e_2 , \dots , e_6 \}$  
around $z$. 
Then, as in \cite[pp.~211]{Tian00}, 
$  \nabla_{\frac{\partial}{\partial r}} \frac{\partial}{\partial r}  = 0
  , 
  \nabla_{\frac{\partial}{\partial r}} X 
 = (\xi' r + \xi ) \phi ( \zeta / r ) 
\frac{\partial}{\partial r} , 
 \nabla_{e_{i}} X
  = \xi r e_{i} \phi \frac{\partial}{\partial r} + \xi \phi \sum_{j=1}^{6} b_{ij} e_{j} 
 \, (2\leq i \leq 6)$, 
$| b_{ij} - \delta_{ij}| \leq C' c(z) r^2$, and 
$C'$ is a positive constant. 
Plugging these into \eqref{fvar}, we obtain 
\begin{equation}
 \begin{split}
  \int_{Z} &|F_{A}^{\perp}|^2  (  \xi' r + 2 \xi  + C' c(z) r^2 \xi ) \phi \,dV_{g} 
  - \int_{Z} |[ u , \bar{u} ]^2 ( \xi' r + 6 \xi  + C' c(z) r^2 \xi) \phi \, dV_{g}  \\
 & \qquad = 4 \int_{Z} \left\{ \xi' r \phi \left| 
\frac{\partial}{\partial r} \lfloor F_{A}^{\perp}  \right|^2 
+ \xi r \left\langle \frac{\partial}{\partial r} \lfloor  F_{A}^{\perp} , 
\nabla \phi  \lfloor F_{A}^{\perp} \right\rangle  \right\}  \\ 
 & \quad \qquad 
   -12  \int_{Z} \left\{ \xi' r \phi \left| 
\frac{\partial}{\partial r} \lfloor  [u, \bar{u}] \right|^2  
+ \xi r \left\langle 
\frac{\partial}{\partial r} \lfloor  [u, \bar{u}] , 
\nabla \phi  \lfloor [u, \bar{u}]
 \right\rangle  \right\}  .
 \end{split}
\label{m2}
\end{equation}

Let $\chi (r)$ be a function which is smooth and satisfies 
$\chi (r) =1 $ for $r
\in [0 ,1] $, $\chi (r) = 0 $ for $r \in [ 1 + \varepsilon , \infty )$ 
where $ \varepsilon >0 $,  and $\chi' (r) \leq 0$.   
We choose $\xi (r) = \xi_{\tau} (r) := \chi ( r / \tau)$ for $\tau 
\in [\sigma , \rho]$. 
Then we obtain 
$ \tau \frac{\partial}{\partial \tau} ( \xi_{\tau} (r)) = - r \xi_{\tau}' (r)$. 
From this with \eqref{m2}, we get  
\begin{equation*}
 \begin{split}
&\frac{\partial}{\partial \tau} 
 \left( \tau^{-2} e^{a \tau^2} \int_{Z} \xi_{\tau} \phi | F_{A}^{\perp}|^2 \right) 
- \tau^{4}  \frac{\partial}{\partial \tau} 
\left( \tau^{-6} e^{a \tau^2} \int_{Z} \xi_{\tau} \phi | [u, \bar{u}]|^2
  \right) \\ 
& \quad \geq 4 \tau^{-2} e^{a \tau^2} \frac{\partial }{\partial \tau } 
 \left( \int_{Z} \xi_{\tau} \phi \left| 
 \frac{\partial}{\partial r} \lfloor F_{A}^{\perp}  \right|^2 \right)  
 -12 \tau^{-2} e^{a \tau^2} \frac{\partial }{\partial \tau } 
 \left( \int_{Z} \xi_{\tau} \phi \left| \frac{\partial}{\partial r} 
   \lfloor  [u, \bar{u}]  \right|^2 \right)  \\ 
 &  \qquad \qquad 
-4 \tau^{-3} e^{a \tau} \int_{Z} \xi_{t} r 
 \left\langle 
\frac{\partial}{\partial r} \lfloor  F_{A}^{\perp} , 
\nabla \phi  \lfloor F_{A}^{\perp}
 \right\rangle \\
& \qquad  \qquad  \qquad
  +12  \tau^{-3} e^{a \tau} \int_{Z} \xi_{t} r 
 \left\langle 
\frac{\partial}{\partial r} \lfloor  [u, \bar{u}] , 
\nabla \phi  \lfloor [u, \bar{u}] 
 \right\rangle ,\\ 
 \end{split}
\end{equation*}
where $a$ is a positive constant with $a \geq C' c (z)$. 
Then, 
integrating this on $\tau$ and letting $\varepsilon$ go to zero, 
we obtain 
\begin{equation*}
 \begin{split}
  & \rho^{-2} e^{a \rho^2} \int_{B_{\rho} (z)} \phi \, 
   m (A ,u) dV_{g} 
 -  \sigma^{-2} e^{a \rho^2} \int_{B_{\sigma} (z)} \phi \,  
   m (A ,u) dV_{g} \\ 
  & \quad \geq \int_{B_{\rho} (z) \setminus B_{\sigma} (z)} r^{-2} e^{a
  r^2} \phi \left(  4 \left| \frac{\partial}{\partial r} \lfloor F_{A}^{\perp} 
  \right|^2 - 12  \left| \frac{\partial}{\partial r} \lfloor  [u,
  \bar{u}] \right|^2 \right)  dV_{g}  \\
  & \qquad  - 4 \int_{\sigma}^{\rho}  \tau^{-3} e^{a \tau} \int_{B_{\tau}
  (z)} r  \left| \frac{\partial}{\partial r} \lfloor  F_{A}^{\perp}
  \right|  \, \left| \nabla \phi \lfloor F_{A}^{\perp} \right| 
   dV_{g} d \tau \\
  & \qquad  \quad  - 12 \int_{\sigma}^{\rho}  \tau^{-3} e^{a \tau} \int_{B_{\tau}
  (z)} r \left|  \frac{\partial}{\partial r} 
   \lfloor  [u, \bar{u}]  \right|  \, 
\left| \nabla \phi \lfloor  [u, \bar{u}]  \right| 
   dV_{g} d \tau \\
& \qquad \quad \quad + \int_{\sigma}^{\rho} 8 \tau^{-3} e^{a \tau^2}
  \int_{B_{\tau} (z)} \phi | [ u, \bar{u}]|^2 dV_{g} d \tau \\ 
 \end{split}
\end{equation*} 
Then taking $\phi \equiv 1$ gives Proposition \ref{th:monotone}. 
\end{proof}

\section{Estimates}
\label{sec:est}

\subsection{A bound on $\bf u$}

In this subsection, we derive a bound on $u$.  
We use ideas by Mares  in the study of the Vafa--Witten
equations in his Ph.D thesis \cite{Ma}.  
\begin{proposition} 
Let $(A,u)$ be a D--T instanton 
of 
a unitary vector bundle $E$ over a compact K\"{a}hler threefold $Z$. 
Then we have 
\begin{equation}
 || u ||_{L^{\infty}} 
\leq C || u ||_{L^2} , 
\end{equation}
where $C >0 $ is a positive constant which depends only on $Z$. 
\label{estu}
\end{proposition}

\begin{proof}
This is basically a rephrasing of \cite[Th.~3.1.1]{Ma}. 
First, we use the following form of the Weitzenb\"{o}ck formula 
(see \cite[Th.~8.17]{LM}): 
\begin{equation}
\bar{\partial}_{A} \bar{\partial}_{A}^{*} 
u =  \nabla_{A}^{*} \nabla_{A} u + s u + [i \Lambda F^{1,1}_{A} , u], 
\label{weit}
\end{equation}
where $s$ is the scalar curvature of the metric $g$, 
and $\Lambda = ( \wedge \omega)^{*}$. 
Then, by using Eqs. \eqref{DT1} and \eqref{DT2}, we get 
$\langle u , \nabla_{A}^{*} \nabla_{A} u \rangle 
= - s | u |^2  - |[u, \bar{u}]|^2$. 
Hence we obtain $\langle u , \nabla_{A}^{*} \nabla_{A} u \rangle 
\leq C | u |^2$. Then we invoke 
\cite[Th.~3.1.2]{Ma} to deduce that $|| u ||_{L^{\infty}} 
\leq C || u ||_{L^2}$. 
\end{proof}

From Proposition 3.1, we immediately get the following.  
\begin{Corollary}
Let $(A,u)$ be a D--T instanton 
of 
a unitary vector bundle $E$ over a compact K\"{a}hler threefold $Z$.  
Assume a bound on the $L^2$-norm of
 $u$. 
Then for any $z \in Z$ and $\varepsilon > 0$, 
there exists a number $r_0 >0$ such that
 for all $0< r \leq r_0$ we have 
$\frac{1}{r^2} \int_{B_{r} (z)} | [ u, \bar{u}] |^2 \,d V_{g} \leq \varepsilon$. 
\label{coru}
\end{Corollary}

\subsection{Curvature estimate}

In this subsection, we prove an $\varepsilon$-regularity theorem for the
Donaldson--Thomas instantons on compact K\"{a}hler threefolds.  

\begin{proposition}
Let $(A,u)$ be a D--T instanton 
of 
a unitary vector bundle $E$ over a compact K\"{a}hler threefold $Z$.  
Then there exist constants $\varepsilon >0$ 
and $C > 0 $
 which depend only on $Z$ such that 
for any $z \in Z$ and $0 < r < r_z$, where $r_z$ is the constant in
 Proposition \ref{th:monotone},  if 
$ \frac{1}{r^{2}} \int_{B_{r}(z)} |F_{A}|^2 \, dV_g \leq \varepsilon$ 
and $\int_{B_{r}(z)} |u|^2 dV_{g} < \varepsilon$
then 
$$ | F_{A} | (z) 
\leq \frac{C}{r^2} \left( \frac{1}{r^{2}} 
\int_{B_{r}(z)} |F_{A}|^2 \,  dV_g  \right)^{\frac{1}{2}} 
+ \frac{C \varepsilon}{r^2}. 
$$
\label{th:est}
\end{proposition}

\begin{proof}
The proof also goes almost identically to that of
 \cite[Th.~2.2.1]{Tian00} except that we have the extra
 nonlinear term coming from the Higgs fields. 
\begin{lemma}
Let $(A,u)$ be a D--T instanton 
of 
a unitary vector bundle $E$ over a compact K\"{a}hler
 threefold $Z$.  
Let $z \in Z$ and $0 < r < r_z$. 
Suppose that $\int_{B_{r}(z)} |u|^2 \, d V_{g} \leq \varepsilon$. 
Then the following holds: 
\begin{equation}
 \Delta  |F_{A}| \geq - C \varepsilon  
 - C' |F_{A}| 
  - C'' |F_{A} |^2,
\label{eq:ineqLap}
\end{equation}
where $C, C', C'' >0 $ are constants which depend only on $Z$. 
\end{lemma}

\begin{proof}
We use the following form of the Weitzenb\"{o}ck formula: 
\begin{equation*}
 \nabla_{A}^{*} \nabla_{A} \varphi 
 = \Delta_{A} \varphi 
 + R(g) \# \varphi + F_{A} * \varphi ,  
\end{equation*}
where $\varphi \in \Omega^{p} ( Z , \mathfrak{u} (E))$, $R(g)$ is the
 Riemannian curvature of $g$, and $\#$ and $*$ are multi-linear maps 
(see \cite[pp.~214]{Tian00} for explicit expressions). 
Then we get 
\begin{equation*}
 \begin{split}
 \Delta |F_{A}|^2 
 &= 2 | \nabla_{A} F_{A} |^2 - 2 \langle \nabla_{A}^{*} 
  \nabla_{A} F_{A} , F_{A}
  \rangle \\ 
 &=  2 | \nabla_{A} F_{A} |^2  
- 2 \langle D_{A} D_{A}^{*} F_{A} 
   + R (g) \#   F_{A} +  F_{A} * F_{A} , F_{A} 
  \rangle  \\ 
 &=  2 | \nabla_{A} F_{A} |^2 - 2 \langle \Lambda^2 
  D_{A} D_{A}^{*} [ u , \bar{u}] 
   + R (g) \#   F_{A} +  F_{A} * F_{A} , F_{A} 
  \rangle  \\  
  &\geq 2 | \nabla_{A} F_{A} |^2   
     - C || [u , \bar{u}] ||_{L^{\infty}} \, | F_{A} | 
      -  C' |F_{A}|^2  
   -  C'' |F_{A} |^3 , \\
  &\geq 2 | \nabla_{A} F_{A} |^2   
   -  C \varepsilon |F_{A}| 
    -  C' |F_{A}|^2  -  C'' |F_{A} |^3 , \\
 \end{split}
\end{equation*}
where we used the equation 
$(1 + *\wedge \omega)F_{A} = \Lambda^2 [u , \bar{u}]$, 
the Bianchi identity $D_{A} F_{A} =0$,  
and Proposition \ref{estu}. 
Thus, we get 
$ \Delta |F_{A}| \geq 
- C \varepsilon  -  C' |F_{A}| - C'' |F_{A} |^2 $.  
\end{proof}

Next, we put 
$ f (\rho) := ( r - 2 \rho )^2 \sup_{x \in B_{\rho} (z)} |F_{A}| (x)$, 
where $\rho \in [ 0, r /2]$. 
This function is continuous; thus it attains its
 maximum at some $\rho_{0} \in [ 0, r /2]$.

\begin{lemma} 
Let $z \in Z$ and $0 < r < r_z$. 
Suppose that $\frac{1}{r^2} \int_{B_r (z)} |F_{A}|^2 dV_{g} \leq 
\varepsilon$ and $\int_{B_{r} (z)} |u|^2 \, dV_{g} \leq \varepsilon $ 
for $\varepsilon  > 0$ sufficiently small. 
Then 
$f (\rho_{0}) \leq 64 $.  
\label{lem:f}
\end{lemma}

\begin{proof}
We put 
$ b = \sup_{x \in B_{\rho_{0}} (z)} |F_{A}| 
(x) = |F_{A}| (x_{0})$, 
and take $\sigma = ( r - 2 \rho_0) /4$. 
Then, we get 
\begin{equation*}
 \begin{split}
   \sup_{x \in B_{\sigma} (x_{0}) }
   & |F_{A}| (x) 
   \leq \sup_{x \in B_{\rho_0 + \sigma} (z)} |F_{A}| (x)\\
  &\leq \frac{( r - 2\rho_{0})^2}{( r - 2 \rho_0 - 2 \sigma )^2} 
    \sup_{x \in B_{\rho_0 } (z)} |F_{A}| (x) = 4 b.  
 \end{split}
\end{equation*}
We now suppose for a contradiction that $f (\rho_{0}) > 64$. 
Clearly we have $\sigma \sqrt{b} \geq 2$. 
We then take $\ell = \text{max} \{ b , 1 / r_{0}^{2} \}$, 
where $r_0$ is the constant in Corollary \ref{coru}.  
We define a new metric by $\tilde{g} := \ell g$, and rescale $(A,u)$ as  
$(A, \tilde{u}) = (A , \ell^2 u)$.  
Then $(A, \tilde{u})$ is a D--T instanton with respect to the metric 
$\tilde{g}$, 
and we get 
$ |F_{A}|_{\tilde{g}} 
= \ell^{-1} |F_{A}|$, where $|F_{A}|^{2}_{\tilde{g}}$
 is the energy density with respect to $\tilde{g}$. 
We then obtain  
$ \sup_{x \in B_{2}(x_{0} ,\tilde{g})} |F_{A}|_{\tilde{g}} 
(x) \leq 4 $.

On the other hand, 
from \eqref{eq:ineqLap}, we have 
$ \Delta_{\tilde{g}} |F_{A}|_{\tilde{g}} \geq 
 -C \varepsilon - C' 
  |F_{A}|_{\tilde{g}} - C'' |F_{A}|_{\tilde{g}}^2 $. 
Thus, we get $
 \Delta_{\tilde{g}} |F_{A}|_{\tilde{g}} \geq 
 - C \varepsilon - C' 
  |F_{A}|_{\tilde{g}}$ on $B_2 (x_0 , \tilde{g})$. 
Hence, by the mean value theorem 
(see e.g. \cite[Th.~9.20]{GT}), we obtain 
\begin{equation*}
 1 + C \varepsilon  
= |F_{A}|_{\tilde{g}} ( x_0) 
 + C \varepsilon 
\leq C 
\left( \int_{B_{1}( x_{0} , \tilde{g})} 
|F_{A}|^{2}_{\tilde{g}} \, dV_{\tilde{g}} \right)^{\frac{1}{2}} 
 + C' \varepsilon . 
\end{equation*}
Furthermore, from Proposition \ref{th:monotone} and Corollary \ref{coru}, we get 
\begin{equation*}
 \begin{split}
 &\int_{B_{1} ( x_{0} , \tilde{g})} |F_{A}|^{2}_{\tilde{g}} \, dV_{\tilde{g}} 
   = ( \sqrt{\ell})^{2} \int_{B_{\frac{1}{ \sqrt{\ell}}} 
   ( x_{0} , g) } |F_{A}|^2 \,  dV_{g} \\ 
  & \quad = ( \sqrt{\ell})^{2} \int_{B_{\frac{1}{ \sqrt{\ell}}} 
   ( x_{0} , g) } m(A,u) dV_{g} 
   + 3 ( \sqrt{\ell})^{2} \int_{B_{\frac{1}{ \sqrt{\ell}}} 
   ( x_{0} , g) } |[ u , \bar{u} ]|^2 \,  dV_{g} \\ 
  & \quad \leq \frac{1}{(r/2)^2} e^{a (r/2)^2} 
  \int_{B_{\frac{r}{2} } ( x_{0} , g) } 
  m(A,u)  dV_{g} 
   + 3 ( \sqrt{\ell})^{2}  \int_{B_{\frac{1}{ \sqrt{\ell}}} 
   ( x_{0} , g) } |[ u , \bar{u} ]|^2 \,  dV_{g} \\  
  & \quad \leq 2^2 e^{(a /4) r^2 } \varepsilon 
 + 3 \varepsilon .  
 \end{split}
\end{equation*} 
Hence, we obtain  
$ 1 \leq C \left( 2^{2} e^{(a/4) r^2} \varepsilon 
 + 3 \varepsilon \right)^{\frac{1}{2}} + C' 
 \varepsilon $, 
but this contradicts the assumption that 
$\varepsilon$ is sufficiently small.  
Thus Lemma \ref{lem:f} holds. 
\end{proof}

From Lemma \ref{lem:f}, we obtain 
$\sup_{x \in B_{r /4 }(z)} r^2 |F_{A}| 
 \leq 4 f (\rho_0) \leq C $. 
We again define a new metric by $g' := r^{-2} g$, 
and rescale $(A,u)$ as $(A , u') = (A, r^{-2} u)$. 
Then we get 
$|F_{A}|_{g'} \leq C$. 
Combining this with \eqref{eq:ineqLap}, 
we obtain  
$ \Delta_{g'} |F_{A}|_{g'} \geq 
 -C \varepsilon  
 - C'  |F_{A}|_{g'} $. 
Hence the mean value theorem again implies Proposition \ref{th:est}. 
\end{proof}

From Proposition \ref{th:est} and a result by Uhlenbeck
\cite[Th.~2.7]{Uh1}, we get the following. 
\begin{Corollary}
There exist constants $\varepsilon >0$, $C > 0$, 
and $r_{\varepsilon} >0$ such that 
for any $z \in Z$ and  $0< r < r_{\varepsilon}$, if $(A,u)$
 is a D--T instanton over $B_{r} (z)$ with 
$r^{-2} \int_{B_{r} (z)} |F_{A}|^2 dV_{g} \leq \varepsilon$, 
then there exists a gauge transformation $\sigma$ over $B_{r} (z)$ 
such that 
$d^{*} \sigma (A) =0$ and 
$|| \sigma (A) ||_{L^{\infty} (B_{r} (z))} \leq C 
||F_{\sigma (A)} ||_{L^{\infty} (B_{r} (z))}$. 
\label{cuh}
\end{Corollary}

\section{A weak convergence. }
\label{sec:weak}

In this section, using the results in Sections \ref{sec:mono} 
and \ref{sec:est}, we 
prove Theorem \ref{th:main}. 
We basically 
follow the proof of a previous result for Yang--Mills connections by Nakajima
 \cite{Nakajima88} 
(see also \cite[Prop.~3.1.2]{Tian00}).

First, we take $\varepsilon$ as in Corollary \ref{cuh}, and consider a set 
\begin{equation*}
S =  \bigcap_{\delta > r > 0} 
\{ z \in Z \, : \,  \liminf_{n \to \infty} 
\frac{1}{r^2} \int_{B_{r}(z)} 
\{ |F_{A_n}^{\perp}|^2  - | [u_{n} , \bar{u}_{n} ]|^2 \}
\, dV_g \geq \varepsilon  \}, 
\end{equation*}
where $\delta$ is the injective radius of $(Z,g)$. 
One can easily check that this $S$ is closed.

\begin{lemma} 
The real two-dimensional Hausdorff measure of $S$ is finite. 
\end{lemma}

\begin{proof}
Let $K$ be a compact subset of $Z$. 
We take a covering $\{ B_{2 \delta} (z_\alpha ) \} \, (\alpha = 1, \dots , N)$ 
of $S \cap K$, where 
$ z_{\alpha} \in S \cap K$, and 
$B_{ \delta} (z_\alpha ) \cap  
 B_{ \delta} (z_\beta ) = \emptyset $ for $\alpha \neq \beta$.  
Then for sufficiently large $n$, $r^{-2} 
\int_{B_{r} (z_{\alpha})} m (A_{n} , u_{n} ) dV_{g} \geq \frac{\varepsilon}{2}$ 
for $\alpha = 1 , \dots , N$. 
Thus, we obtain 
\begin{equation*}
 \sum_{\alpha} r^2 
 \leq \frac{2}{\varepsilon} \sum_{\alpha} 
 \int_{B_{\delta} (z_{\alpha})} m (A_{n} ,u_{n}) dV_{g} 
 \leq \frac{2}{\varepsilon} 
 \int_{Z} m (A_{n} ,u_{n}) dV_{g} 
 \leq \frac{C}{\varepsilon}. 
\end{equation*}
Hence, the real two-dimensional Hausdorff measure of $K \cap S$ 
is finite. 
\end{proof}

We next take a point $z \in Z \setminus S$. 
By the definition of the set $S$, we can find a number 
$N \in \N$ and a radius $r' > 0 $ 
such that  
$ \frac{1}{{r}^2} \int_{B_{r} (z)} m(A_{n} , u_{n}) dV \leq 
  \varepsilon 
$ for any $0< r < r'$ and $n \geq N$.   
We take $r \leq  \text{min} \, \{ r' , r_z , r_0 , r_{\varepsilon} 
\}$. 
From Corollary \ref{coru}, 
we have $\frac{1}{r^2} \int_{B_{r}(z)} 
|[ u_n , \bar{u}_n]|^2 d V_{g} \leq \varepsilon $ 
 for all $n \in \N$. 
Hence we get 
$$ 
\frac{1}{r^2} \int_{B_{r}(z)} 
|F_{A_n}|^2 d V 
\leq 
\frac{2}{r^2} \int_{B_{r}(z)} 
|[ u_n , \bar{u}_n]|^2 d V_{g} + \varepsilon 
\leq 3 \varepsilon . $$
Therefore, from Corollary \ref{cuh}, 
there exists a Coulomb gauge  
$\sigma_{n}$  such
 that  
$ d^{*} \sigma_{n} (A_{n}) = 0 $ on $B_{r} (z)$ 
with  
$ || \sigma_{n} ( A_{n}) ||_{L^{\infty} (B_{r} (z))} 
 \leq C || F_{\sigma_{n} ( A_{n})} ||_{L^{\infty} (B_{r} (z) )}
$.
Since Eqs. \eqref{DT1} and  \eqref{DT2} are gauge invariant, 
each 
$ ( \sigma_{n} (A_{n}) , \sigma_{n} (u_{n})) $ satisfies the
 Donaldson--Thomas equations. 
Furthermore, Eqs. \eqref{DT1} and  \eqref{DT2} with 
$ d^{*} \sigma (A) = 0 $ form an elliptic system; 
thus,  by standard elliptic theory, 
we also get uniform bounds on the derivatives of $(\sigma_{n} (A_{n}) , 
\sigma_{n} (u_{n}))$. 
Hence, there exists a subsequence which converges to a 
 D--T instanton on $B_{r/2} (z)$ in smooth
 topology.

We then patch the gauges above together 
by using arguments in \cite[\S 4.4.2]{DK} (see also 
\cite[\S 3]{Uh2} )
to get 
a sequence of gauge transformations $\sigma_{j}$ on $Z \setminus S$ 
and  a subsequence $\{ (A_{n_j} , u_{n_{j}}) \}$ 
 such
 that $\sigma_{j} ( A_{n_{j}} , u_{n_{j}})$ converges to a D--T instanton
 on $Z \setminus S$. 
We omit the detail here, since it becomes a formal repetition. 
\qed


\section{A convergence}
\label{sec:conv}

In this section, 
we prove Theorem \ref{th:conv}, which can be thought of as 
an $n/2$-compactness theorem for the Donaldson--Thomas instantons 
on compact K\"{a}hler threefolds. 
This sort of analysis 
for the Yang--Mills and the coupled Yang--Mills fields 
were studied by Sibner \cite{Sibner85}, 
and the convergence results for them were obtained by Zhang 
 \cite{MR1926448}, \cite{MR2086545}.  
The proof goes by analysing the singular set of a limit D--T instanton
in the same way as in \cite{MR1926448}, \cite{MR2086545}.

\paragraph{Convergence of measures and the structure of singular sets.}

First, from Proposition \ref{th:est} 
and the H\"{o}lder inequality, 
we immediately obtain the following.  
\begin{Corollary}
Let $(A,u)$ be a D--T instanton 
of 
a unitary vector bundle $E$ over a compact K\"{a}hler threefold $Z$.  
Then there exist constants $\varepsilon > 0$ 
and $C>0$ such that 
for any $z \in Z$ and $0 < r < r_z$, if 
$\int_{B_{r}(z)} |F_{A}|^3 \, dV_g \leq \varepsilon $ 
and $\int_{B_{r}(z)} |u|^2 \, dV_{g} \leq \varepsilon$,  
then 
$$|F_{A}| (z) \leq \frac{C}{r^2} \left(
\int_{B_{r}(z)} |F_{A}|^3 \,  dV_g
 \right)^{\frac{1}{3}} 
 + \frac{C \varepsilon}{r^2}
. $$ 
\label{cor_curv31}
\end{Corollary}

We then prove the following reproduction of \cite[Th.~4.2]{MR2086545} 
(see also \cite[Th.~4.2.3]{Tian00}). 
\begin{proposition}
 Let $\{ (A_n , u_n) \}$ be a sequence of D--T instantons of
 a unitary vector bundle $E$ over a compact K\"{a}hler threefold $Z$. 
We assume that $ \int_{Z} |F_{A_n}|^3 dV_{g} \leq C$, 
and that  $\int_{Z} | u_{n} |^2 dV_{g}$ are also uniformly
 bounded. 
Then, there exist a subsequence $\{ (A_{k} , u_{k}) \}$, a
 sequence of gauge
 transformations $\{ \sigma_{k} \}$,  and a finite set of points
 $T = \{ z_{\alpha} \}_{\alpha=1}^{\ell} \subset Z$ such that
 $\sigma_{k} (A_{k} , u_{k})$
 converges to a D--T instanton $(A,u)$ over  $Z \setminus T$.

Moreover, for each $\alpha= 1, \dots , \ell$ there exists a positive
 constant 
 $\theta_{\alpha} > 0$ such that 
$$ |F_{A_{k}}|^3  d V_{g} \to 
 |F_{A}|^3  d V_{g} 
+ \sum_{\alpha=1}^{\ell} \theta_{\alpha} \delta_{z_{\alpha}} $$
weakly in the sense of the Radon measure, 
where $\delta_{z_{\alpha}}$ is the
 Dirac measure at $z_{\alpha}$.  
\label{convRadon}
\end{proposition}

\begin{proof}
We follow the proof of \cite[Th.~4.2]{MR2086545}.
First we put 
\begin{equation*}
  T:= \bigcap_{\delta > r > 0}
\{ z \in Z \, | \, 
 \liminf_{n \to \infty} 
\int_{B_{r} (z)} |F_{A_n}|^3 dV_{g} \geq \varepsilon \}.
\end{equation*}
One can easily check that the set $T$ is closed, and 
can prove the following two lemmas in the same way as in Theorem 1.1.

\begin{lemma}
The zero-dimensional Hausdorff measure of $T$ is finite. 
\end{lemma}

\begin{lemma}
There exist a subsequence $\{ (A_{k} , u_{k} )\}$ of 
$\{ (A_{n} , u_{n}) \}$ and a sequence of gauge transformations $\{ \sigma_{k} \}$ such that 
$\sigma_{k} (A_{k} ,u_{k})$ converges to $(A , u)$ outside $Z \setminus T$. 
\end{lemma}

We then consider the Radon measures 
$ \mu_{k} := |F_{A_k}|^3 dV_{g}$. 
By taking a subsequence if necessary, 
$\mu_k$ weakly converges to a Radon measure $\mu$ on $Z$. 
We write $ \mu = |F_{A}|^3 dV_{g} + \nu $, 
where $\nu$ is a non-negative Radon measure on $Z$. 
Since  the support of $\nu$ is in  $T$, 
we write $\nu = \sum \theta_{\alpha} \delta_{y_{\alpha}}$, where
 $\theta_{\alpha} \geq 0$. 
In fact, $\theta_{\alpha}$ is positive, because, 
by using a cut-off function  $\chi \in C^{\infty} (Z)$ 
with $\chi (z) = 1$ on $B_{r} (z_{\alpha}) $ and
 $\chi (z) = 0$ on $Z \setminus B_{2r} (z_{\alpha})$, 
where $B_{2r} (z_{\alpha})$ is 
a geodesic ball of radius $2r$ 
with centre at
 $z_{\alpha}$ so that 
 $T \cap B_{2r} (z_{\alpha}) = \{ z_{\alpha} \}$, 
we obtain 
\begin{equation*}
\begin{split}
  \varepsilon 
 &\leq \liminf_{k \to \infty} \int_{B_{r}(z_{\alpha})} 
  | F_{A_k} |^3 dV_{g} \\
 &\leq \lim_{k \to \infty} \int_{Z} 
 \, \chi | F_{A_{k}} |^3 dV_{g} 
  \leq \theta_{\alpha} 
 + \int_{B_{2r}(z_{\alpha})} 
 |F_{A}|^3  dV_{g} . 
 \end{split}
\end{equation*}
Thus, by taking $r \to 0$, we get 
$\theta_{\alpha} \geq \varepsilon >0$. 
\end{proof}

We next 
take normal coordinates $(\zeta_1 , \dots , \zeta_6)$ around
$z_{\alpha}$, 
and denote by $B (\zeta, \rho)$ an open ball 
of radius $\rho$ with centre at $\zeta$ 
in the normal coordinates. 
Imitating \cite[(4.9)]{MR2086545}, 
we define a function 
$$\tilde{E} ( k , \rho) 
:= \sup_{\zeta \in B (0,r)} \int_{\exp_{z_{\alpha}} (B ( \zeta , \rho))} 
 |F_{A_k}|^3  dV_{g} $$ 
for  $0 \leq \rho \leq r$. 
The function $\tilde{E} ( k , \rho)$ is continuous and
 non-decreasing in $\rho$ and $\tilde{E} ( k , 0) = 0$.

From the definition of $T$, 
we have 
$ \tilde{E} ( k , r) \geq 
 \int_{B_{r} (z_{\alpha})} |F_{A_k}|^3 dV_{g} 
 \geq \frac{3 \varepsilon}{4} $ 
for $k$ sufficiently large. 
Since $\tilde{E} ( k , r)$ is continuous, 
there exist $0 < \rho_{k} < r$ and $\zeta_{k} \in \overline{B ( 0, r)}$ such that 
$ \tilde{E} ( k , \rho_{k}) 
=  \int_{\exp_{z_{\alpha}} (B ( \zeta_{k} , \rho_{k}))} 
 | F_{A_{k}} |^3   dV_{g} = \frac{\varepsilon}{2} $. 
Since $T \cap B_{2r} (z_{\alpha}) 
= \{ z_{\alpha }\}$, 
$\rho_{k} \to 0 , \, \eta_{k} \to 0$ as $k \to \infty$.

We now define 
$A_{k, \rho_{k}} := 
\tau^{*}_{\rho_{k}} \exp^{*}_{z_{\alpha}} A_{k} ,\,\, 
u_{k ,  \rho_{k}} := \rho_{k}^{-2} \tau^{*}_{\rho_{k}} \exp^{*}_{z_{\alpha}} u_{k}$,  
where
 $\tau_{\rho_{k}} (v) := \zeta_{k} + \rho_{k} v$ for $v \in T_{z_{\alpha}}Z$. 
Then, 
$(A_{k,  \rho_{k}} ,  u_{k , \rho_{k}})$ satisfies 
the Donaldson--Thomas equations on 
$T_{z_{\alpha}} Z$ with respect to a  metric 
$g_{k} := \rho_{k}^{-2} \tau_{k}^{*} \exp_{z_{\alpha}}^{*} g $. 
Moreover, we have  
$
 \int_{T_{z_{\alpha}} Z} |F_{A_{k,  \rho_{k}}}|^3  
 dV_{g_{k}} 
 =  \int_{B_{2r} (z_{\alpha})} |F_{A_{k}}|^3
 dV_{g} \leq C $.

The following, a replication of \cite[Th.~4.3]{MR2086545}, also holds for
the Donaldson--Thomas instantons case. 
\begin{proposition}
 Let $(A_{k} , u_{k})$ and $T$ be as in Proposition \ref{convRadon},
 and let $z_{\alpha} \in T$. 
Then there exists  
a subsequence of $\{ (A_{k , \rho_{k}} , u_{k , \rho_{k}})\}$, 
which converges 
to a smooth D--T instanton $(B , v)$ 
of the trivial bundle over $(T_{z_{\alpha}} Z , g_{z_{\alpha}})$ with 
$|F_{B}| \neq 0$ and $\int_{T_{z_{\alpha}} Z} 
|F_{B}|^3 d V_{g_{z_{\alpha}}} \leq \theta_{\alpha}$. 
\label{scale}
\end{proposition}

\begin{proof}
The proof is formally the same as that of \cite[Th.~4.3]{MR2086545}. 
We have 
\begin{equation}
\begin{split}
\tilde{E} (k , &\rho_{k}) 
= \int_{B ( 0, 1)} |F_{A_{k , \rho_k}}|^{3}
 dV_{g_{k}} \\ 
 &= \sup_{\zeta\in \tau_{k}^{-1} (B ( 0, r))} 
 \int_{B (\zeta ,1)} |F_{A_{k , \rho_k}}|^{3}
 dV_{g_{k}} 
= \frac{\varepsilon}{2} . 
\label{ep2}
\end{split}
\end{equation}
Thus, from Corollary \ref{cor_curv31}, 
if $\varepsilon$ is small, 
we get 
$|F_{A_{k ,\rho_{k}}}| 
\leq C \varepsilon^{\frac{1}{3}} + C' \varepsilon$. 
Hence, there exist a subsequence $(A_{k',   \rho_{k'}} , u_{k', 
 \rho_{k'}} )$ and a sequence of 
gauge transformations $\{ \sigma_{k'} \}$ such that 
$\sigma_{k'} (A_{k',   \rho_{k'}} , u_{k', 
 \rho_{k'}} )  $ converges to a D--T instanton $(B ,v)$
 on $(T_{z_{\alpha}} Z , g_{z_{\alpha}}) \cong (\C^3 , g_0)$.  
From  \eqref{ep2}, 
we have 
$ \int_{B(0,1)} |F_{B}|^3 dV_{g_{z_{\alpha}}} = \frac{\varepsilon}{2} $. 
Thus, $|F_{B}| \neq 0$. 
Also by Fatou's lemma, 
\begin{equation*}
 \begin{split}
  \int_{T_{z_{\alpha}}Z} |F_{B}|^3 dV_{g_{z_{\alpha}}} 
  &\leq \liminf_{k' \to \infty} 
  \int_{T_{z_{\alpha}} Z} |F_{A_{k'}}|^3 dV_{g_{k'}} \\
  &\leq \theta_{\alpha} + \int_{B_{2r}(z_{\alpha})} |F_{A}|^3 
     dV_{g} .
 \end{split}
\end{equation*}
Thus, taking $r \to 0$, we obtain 
$ \int_{T_{z_{\alpha}} Z} 
 |F_{B}|^3  dV_{g_{z_{\alpha}}} \leq \theta_{\alpha} $. 
\end{proof}

\paragraph{Proof of Theorem \ref{th:conv}.}

From  Proposition \ref{convRadon}, we can find 
a subsequence $\{ (A_{k} ,u_{k})\}$ of 
$\{ ( A_n , u_n ) \}$ and a sequence of gauge transformations $\{
\sigma_{k} \}$ such that $\sigma_k (A_k , u_k)$ converges to a D--T instanton 
over $Z \setminus T$, where $T$ is a finite set of points. 
If $T \neq \emptyset$, then by using Proposition  \ref{scale}, 
we can construct a D--T instanton $(B,v)$ on $\C^3$ with 
$\int_{\C^3} |F_{B}|^3 dV_{g_{0}} < C$ and $|F_{B}| \neq 0$. 
On the other hand, from the Weitzenb\"{o}ck formula \eqref{weit}, $u$
vanishes on $(\C^3 , g_{0})$, namely, D--T instantons are just Hermitian--Einstein
connections on $(\C^3 , g_{0})$. Thus we get a Hermitian--Einstein connection on 
$( \C^3 , \omega_0)$ 
 with 
$ \int_{\C^3} |F_{A}|^{\frac{3}{2}} dV_{g_{0}} \leq C 
$ and $F_{A} \neq 0$. 
However,  this contradicts the following result by Zhang \cite{MR1926448}.   
\begin{theorem}[\cite{MR1926448} Theorem 3.3]
 If $A$ is a Hermitian--Einstein connection 
over $( \C^3 , \omega_0)$ 
 with 
$ \int_{\C^3} |F_{A}|^{3} dV_{g_{0}} \leq C 
$, 
where $g_{0}$ is the standard metric on $\C^3$, 
then $F_{A} \equiv 0 $.  
\label{cor_curv32}
\end{theorem}
Thus, $T = \emptyset$. This proves Theorem \ref{th:conv}. 
\qed


\begin{flushleft}
E-mail: yuuji.tanaka@math.nagoya-u.ac.jp
\end{flushleft}

\end{document}